\documentclass[lettersize,journal]{IEEEtran}
\usepackage{amsmath,amsfonts}
\usepackage{algorithmic}
\usepackage{algorithm}
\usepackage{array}
\usepackage[caption=false,font=normalsize,labelfont=sf,textfont=sf]{subfig}
\usepackage{textcomp}
\usepackage{stfloats}
\usepackage{url}
\usepackage{verbatim}
\usepackage{graphicx}
\usepackage{cite}
\usepackage{todonotes}
\usepackage{supertabular}
\usepackage{tabularx,booktabs}
\newcolumntype{C}{>{\centering\arraybackslash}X} 
\usepackage{lipsum} 
\begin{document}

\title{Beyond Price-Taker: Multiscale Optimization \\ of Wind and Battery Integrated Energy Systems}


\author{Xinhe Chen, Xian Gao, Darice Guittet, Radhakrishna Tumbalam-Gooty, Bernard Knueven, John D. Siirola, David C. Miller and Alexander W. Dowling








\thanks{
Xinhe Chen, Xian Gao$^*$ and Alexander W. Dowling are ($^*$was) with the Department of Chemical and Biomolecular Engineering, University of Notre Dame, Notre Dame, IN 46556, USA. (email: xchen24@nd.edu; email: xgao1@alumni.nd.edu; adowling@nd.edu [corresponding author])

Darice Guittet and Bernard Knueven are with the National Renewable Energy Laboratory, Golden, CO 80401, USA. (email: darice.guittet@nrel.gov; Bernard.Knueven@nrel.gov)

David C. Miller$^*$ was with the National Energy Technology Laboratory (NETL), Pittsburgh, PA 15236, USA. (email: david.miller@olisystems.com). 

Radhakrishna Tumbalam Gooty is an NETL Support Contractor, Pittsburgh, PA 15236 (email: Radhakrishna.Tumbalam-Gooty@netl.doe.gov).

John D. Siirola is with Sandia National Laboratories, Albuquerque, NM 87185, USA. (email: jdsiiro@sandia.gov)

}}

\markboth{Journal of \LaTeX\ Class Files,~Vol.~14, No.~8, August~2021}%
{Shell \MakeLowercase{\textit{et al.}}: A Sample Article Using IEEEtran.cls for IEEE Journals}


\maketitle

\begin{abstract}
Decarbonizing global energy systems requires extensive integration of renewable energy into the electric grid. However, the intermittency and variable nature of wind and other non-dispatchable renewable resources make integration a great challenge. Hybridizing renewables with energy storage to form integrated energy systems (IESs) helps mitigate these concerns by improving reliability and resiliency. This paper systematically studies the limitations of the prevailing price-taker assumption for techno-economic analysis (TEA) and optimization of hybrid energy systems. As an illustrative case study, we retrofit an existing wind farm in the RTS-GMLC test system (which loosely mimics the Southwest U.S.) with battery energy storage to form an IES. We show that the standard price-taker model overestimates the electricity revenue and the net present value (NPV) of the IES up to 178\% and 50\%, respectively, compared to our more rigorous multiscale optimization. These differences arise because introducing storage creates a more flexible resource that impacts the larger wholesale electricity market. Moreover, this work highlights the impact of the IES has on the market via various strategic bidding, and underscores the importance of moving beyond price-taker for optimal storage sizing and TEA of IESs. We conclude by discussing opportunities to generalize the proposed framework to other IESs, and highlight emerging research questions regarding the complex interactions between IESs and markets.

\end{abstract}

\begin{IEEEkeywords}
renewable energy, battery energy storage, production cost models, energy markets, net present value.
\end{IEEEkeywords}

\section*{nomenclature}
\begin{tabular}{p{0.15\textwidth} p{0.30\textwidth}}
\textbf{\textit{Variable}} & \textbf{\textit{Description}}\\
\hline
$\bar{P}^{b}$ & Maximum battery power, MW\\
$\bar{H}$ & Battery capacity, hr\\
$\bar{S}$ & Maximum battery state of the charge (SoC), MWh\\
$p^{w}_{t,i}$ & Wind power generated at time $t \in \mathcal{T},$ and in scenario $ ~i \in \mathcal{I}$, MW\\
$p^{s}_{t,i}$ & Power output from the wind farm to the grid at time $t \in \mathcal{T},$ and in scenario $ ~i \in \mathcal{I}$, MW\\
$p^{c}_{t,i}$ & Power charging the battery at time $t \in \mathcal{T},$ and in scenario $ ~i \in \mathcal{I}$, MW\\
$p^{d}_{t,i}$ & Power charged from the battery at time $t \in \mathcal{T},$ and in scenario $ ~i \in \mathcal{I}$, MW\\
$p_{t,i}$ & Total power sold to the grid at time $t \in \mathcal{T},$ and in scenario $ ~i \in \mathcal{I}$, MW\\
$S_{t,i}$ & Battery SoC at time $t \in \mathcal{T},$ and in scenario $ ~i \in \mathcal{I}$, MWh\\
$E_{t,i}$ & Battery energy accumulated energy throughput at time t $t \in \mathcal{T},$ and in scenario $ ~i \in \mathcal{I}$, MWh\\
$C_{OM}$ & Operation and maintenance (O\&M) cost of the IES, \$\\
$R^{RT}_{t,i}$ & Real-time electricity revenue at time $t \in \mathcal{T},$ and in scenario $ ~i \in \mathcal{I}$, \$\\
$C_{capex}$ & Capital cost of the battery, \$\\
$p^{RT}_{t,i}$ & Total real-time power sold by IES at time $t \in \mathcal{T},$ and in scenario $ ~i \in \mathcal{I}$, MW\\
\quad & \quad \\ 
\textbf{\small \textit{Parameter}} & \small \textbf{\textit{Description}}\\
\hline
$\bar{P}^{w} = 847$ & Maximum wind farm power, MW\\
$f_t \in [0,1] $ & Capacity factor of the wind farm $t \in \mathcal{T}$, dimensionless\\
$C_{OM}^{w} = 42$ & O\&M cost of per unit wind farm power per year, k\$/MW-yr\\
$C_{OM}^{b} \in [19, 70]$ & O\&M cost of per unit battery SoC per year, k\$/MWh-yr\\
$\eta^{c} = 0.95$ & Charging efficiency, dimensionless\\
$\eta^{d} = 0.95$ & Discharging efficiency, dimensionless\\
$\Delta t = 1$ & time resolution, hr\\
$S_{0,i} = 0$ & Initial SoC, MWh\\
$S_{\mathcal{T},i} = 0$ & Final SoC, at scenario $ ~i \in \mathcal{I}$, MWh\\
\end{tabular}
\begin{tabular}{p{0.15\textwidth} p{0.30\textwidth}}
$\delta = 10^{-4}$ & Battery degradation coefficient, dimensionless\\
$C_{capex}^{b} \in [0.8, 3.1] $ & Capital cost of per unit power, k\$/kW\\
$\epsilon = 10^{-3}$ & Renewable production incentive, \$/MWh\\
$\phi = 15.4$ & NPV factor, dimensionless\\
$N = 30$ & Life cycle of the investment, year\\
$r = 0.05$ & Discount rate, dimensionless\\
\quad & \quad \\ 
\textbf{\small \textit{Set}} & \small \textbf{\textit{Description}}\\
\hline
$i \in \mathcal{I}$ & LMP signals (price-taker)/forecasts (time-variant bidding) for scenario $i$\\
$t \in \mathcal{T}$ & Time steps for timescale $t$\\
\end{tabular}

\section{Introduction}
\label{sec:1}
\IEEEPARstart{G}{LOBAL} renewable energy generation has grown 81\% from 2012 to 2022 \cite{worldrenewable}. The United Nations estimates that renewable electricity could provide 65\% of the world’s total electricity supply by 2030 and decarbonize 90\% of the power sector by 2050 \cite{unrenewable}. 
During the past decade, the total installed capacity of wind power in the U.S. has grown from 61 GW to 154 GW \cite{usdoe}. However, the intermittent nature of non-dispatchable renewables poses great challenges to their integration into the electric grid \cite{vargas2019wind}, as adding more renewables can increase grid instability, intermittency, and price volatility\cite{ertugrul2016battery}. Integrated energy systems (IESs) help improve the reliability, economics, and sustainability of renewable power systems \cite{arent2021multi} by exploiting synergies among various technologies such as fossil fuels \cite{kang2011optimal, wang2019operation, chen2020optimal}, nuclear \cite{bragg2020reimagining}, wind \cite{li2017optimal, li2021two, li2020improving}, solar \cite{jiang2020optimal} hydrogen \cite{ruiming2019multi, fabianek2023techno} and energy storage \cite{sioshansi2021energy} to provide multi-products and services to multiple markets. Hybridizing wind with battery storage is especially attractive because of its high specific energy (energy per unit weight) and density (energy per unit volume) of batteries \cite{koohi2020review}. Moreover, battery storage offers fast responses, improved control, and geographical independence\cite{yang2018battery}. Integrating a renewable energy system with battery storage can turn non-dispatchable renewables into dispatchable power and reduce renewable curtailment\cite{denholm2023moving}. Moreover, wind-battery IESs can arbitrage multiple markets and timescales\cite{dowling2017multi}, thus improving the economics of wind resources \cite{bechlenberg2024renewable}. 
\begin{table*}[b]
\caption{Published case studies of renewable-battery integrated energy system. (DA: day-ahead market; RT: real-time market; FR: frequency regulation market)}
\label{tb_ies}
\begin{tabularx}{\textwidth}{@{} l p{3cm} p{2cm} p{1.5cm} p{2.5cm} p{6.5cm} @{}}
\toprule
& Author  & Location & Market & Model & Main Conclusion\\ 
\midrule
& Dicorato et al. \cite{Dicorato2012planning}       & Southern Italian    &DA, RT& PT & High-capacity storage devices improve market performance and comply with the delivery plan.\\
& Moghaddam et al. \cite{Moghaddam2019bptimal}       & MISO    &DA, RT& PT & Voltage control, battery storage lifetime, and battery storage location will impact the average daily profit.\\
& Bechlenberg et al. \cite{bechlenberg2024renewable}     & Belgium/Nordpool  &DA, RT,  FR&PT with bidding &  Sizing the entire wind+battery IES can improve NPV by 50\% compared to sizing battery storage only.\\ 
& Hou et al. \cite{hou2018cooperation}       & Nordpool   &DA, FR& PT with bidding & Participation in multiple markets can increase the profit and help the grid operator stabilize the grid.\\
& He et al.\cite{he2016cooperation}        & PJM   &FR& PT with bidding & Battery storage increases the wind regulation performance score and economic benefits.\\
& Nguyen et al.\cite{nguyen2012new} & PJM   &RT& PT & A dynamic programming framework to optimize the IES operating schedule can improve the profit.\\
\bottomrule
\end{tabularx}
\end{table*}


Design and operation optimization of IESs, including renewable-battery systems \cite{yang2018battery}, is challenging because different design and operation decisions across multiple timescales span from seconds to decades\cite{GAO2022119017}. The price-taker (PT) approximation is the \emph{de facto} standard for optimal sizing and operation of IESs because of its simplicity. It treats electricity locational marginal prices (LMP) as exogenous inputs to the optimization model and assumes that the IES has negligible impact on the market. For example, Sorourifar et al. \cite{sorourifar2018integrated} used multiscale linear programming to optimally size and schedule replacements for grid-scale NaS batteries, considering historical LMP from the California Independent System Operator (CAISO) and market participation decisions in day-ahead (hours), real-time (minutes), and ancillary service markets (seconds). They emphasized that full multiscale market participation increases the net present value (NPV) 4 to 5 times compared to participation in the day-ahead market alone. Dicorato \cite{Dicorato2012planning} performed a life cycle analysis (decades) of a wind-battery IES while considering hourly LMP signals using the PT approximation. They found that high-capacity storage devices enable the integration of wind power into production planning and dispatch and reduce the yearly penalty for unbalanced power by 80\%. Moghaddam et. al \cite{Moghaddam2019bptimal} used a PT model with the Midcontinent Independent System Operator (MISO) historical data to optimally size battery storage for a wind farm. Their optimization model maximizes the profit and shows that the voltage control, the battery storage lifetime, and the battery storage location would impact the average daily profit. Table \ref{tb_ies} summarizes additional wind-battery IES optimization studies.

However, the PT approximation oversimplifies the complex and dynamic electricity market by assuming that the electric grid is an ``infinite bus'' capable of providing and absorbing any (reasonable) amount of power without impacting market prices \cite{GAO2022119017}. This makes the PT approximation easy to use and computationally inexpensive at the cost of ignoring the impact of IES decisions on market outcomes. However, recent studies probe how this assumption can result in inaccurate and misleading conclusions. For example, Emmanuel et al. \cite{emmanuel2022market} showed that the PT approximation fails to predict the suppression of the market price with the deployment of large-scale battery storage. Frew et al. \cite{frew2023analysis} found that the hybridization of nuclear power plant with hydrogen co-production could increase the LMP by 3.5 \$/MWh, and the PT model overestimates the NPV by 200 million\$ compared to rigorous production cost model (PCM) simulations. Similarly, Teixeira et al. \cite{teixeira2012market} developed a nonlinear model to predict how pumped storage hydro (PSH) influences the market clearing price, considering how market power reduces the value of PSH compared to the PT model.

PCMs mimic the price-setting mechanisms executed by Independent System Operators (ISOs) and Regional Transmission Organizations (RTOs). PCMs schedule generators to minimize the total cost of satisfying the forecasted total demand in the entire region/market while considering operational constraints (e.g., maximum/minimum power and ramping rate), forecasted renewable availability, generator costs, transmission limits, and reliability margins. These computationally expensive unit commitment and economic dispatch optimization problems are solved in a multiscale rolling horizon to set prices for multiple market timescales as well as power generation and transmission schedules \cite{batlle2005strategic}. Thus, PCM optimization problems do not necessarily maximize the revenue for individual generators or IESs. PCM simulations help stakeholders and decision-makers understand grid-level performance from system upgrades, e.g., adding transmission capacity, retiring aged generators, or introducing new IESs \cite{sioshansi2021energy}. Recently, Gao et al. \cite{GAO2022119017} proposed a technology-agnostic multiscale computational framework to integrate dynamic optimization of IESs and grid-level PCMs. They highlighted the nuanced ways generator retrofits or replacements impact market outcomes, underscoring the oversimplifications in the PT approximation. Jalving et al. \cite{jalving2023beyond} extended this framework by using machine-learning surrogate models to predict market outcomes as a function of IES design decisions. These surrogates, which explicitly emulate complex IES-market interactions, were embedded in IES design and operation co-optimization. This analysis highlighted market and design conditions, e.g., high marginal costs, where the benchmarked PT optimization becomes unreliable.

Battery storage provides flexibility to the grid to mitigate uncertainty. However, representing this flexibility is non-trivial for IESs that contain non-dispatchable renewables and energy storage. As illustrated in Fig. \ref{fig:ies model}, energy bids reflect the production capacity and marginal costs of a resource, e.g., generator or IES. In the prevailing fossil-centric generation paradigm, bid curves are the generation costs as a function of power output, parameterized by fuel price. Thus, bid curves inform the market about the true costs of ramping and shutting down generators. Within this paradigm, non-dispatchable renewable generators bid a near-zero marginal cost since they consume no fuel. 
ISOs and RTOs use sophisticated forecasts and ancillary services (e.g., reserves, frequency regulation) to mitigate the uncertainty of renewables. Hybridizing renewables with storage enables these IESs to participate similarly to fossil generators. An open question is how these new IESs should strategically bid into the market as dispatchable generators, as their marginal costs are time-varying due to their dependence on the amount of stored energy and forecasted renewable availability. Recently proposed strategies for IES bidding include approaches based on machine learning \cite{cocchi2018machine}, dynamic programming \cite{zheng2023energy}, and game-theory \cite{kian2005bidding}. However, these studies only benchmark bidding strategies under the PT approximation, i.e., they neglect how changing the bidding strategy impacts market outcomes. In this paper, we argue that these interactions are complicated and significant, which motivates a more rigorous analysis of IES-market interactions.


This paper systematically compares the pervasive PT approximation against rigorous multiscale models \cite{GAO2022119017} for optimal sizing and operation of a battery storage system to retrofit a wind farm, resulting in a dispatchable IES. The novel contributions include:
\begin{itemize}
    \item Quantifying the limitations and biases of the PT approximation. 
    \item Providing guidance on conditions when the PT approximation is most unreliable.
    \item Demonstrating the benefits of IES bidding strategies within the context of the market, e.g., using a PCM, which is often neglected in other studies.
\end{itemize}


\section{Methods}
\label{sec:2}

Our goal is to determine the optimal battery energy storage size and operating strategy to retrofit wind farm ``303\_WIND\_1'' at Bus ``Caesar'' in the RTS-GMLC network into a dispatchable IES that participates in the real-time energy market. The remainder of this section describes the mathematical model, case study, and computational implementation summarized in Fig \ref{fig:ies model}.

\subsection{Wind-Battery IES Model}
\begin{figure}[h]
\centering
\includegraphics[width=3.5in]{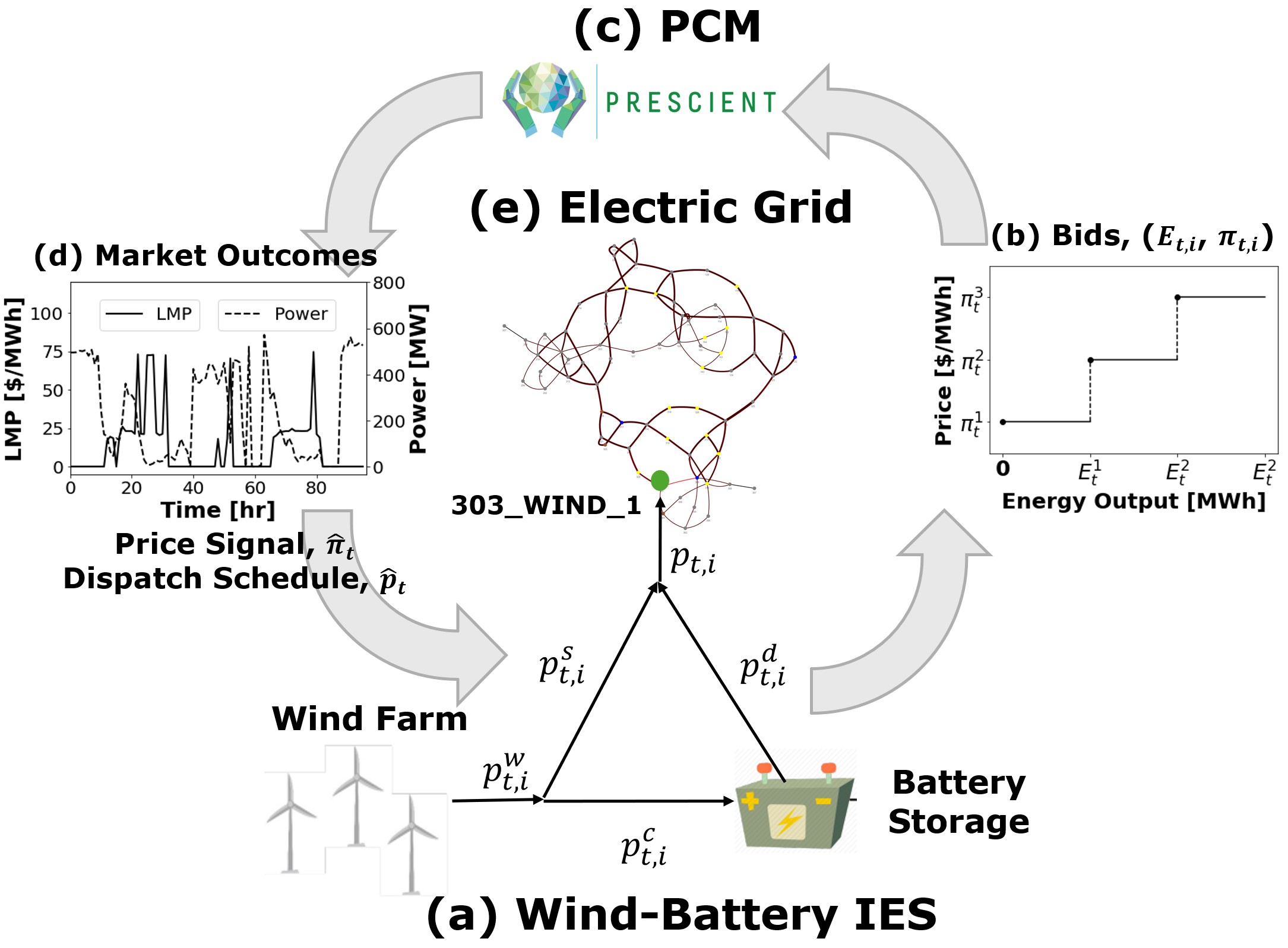}
\caption{At each time step, the wind-battery IES (a) submits its bids (b) to the market. Then the PCM (c) solves unit commitment problems and economic dispatch problems to generate market outcomes (d) for the IES, including the LMP $\hat{\pi}_t$ and dispatch schedule $\hat{p}_{t}$. Then the IES dispatches electricity to the electric grid (e) to satisfy the dispatch schedule.}
\label{fig:ies model}
\end{figure}

Fig. \ref{fig:ies model} summarizes our wind-battery IES system and how it interacts with the electricity market. As the wind farm produces power, the IES can decide to charge the battery storage with power $p^{c}_{t,i}$, and deliver the remaining power $p^{s}_{t,i}$ to the grid. Meanwhile, the IES can also discharge the battery and sell power $p^{d}_{t,i}$ to the grid. The total power delivered to the grid is $p_{t,i}$, which is the sum of $p^{s}_{t,i}$ and $p^{d}_{t,i}$. This work assumes that the IES does not buy power from the grid, because when the LMP is low, the wind farm usually has very high renewable curtailment. 

We define two indices and sets: $t \in \mathcal{T}$ is the set of time points and $i \in \mathcal{I}$ is the set of LMP scenarios/forecasts. The wind power generation $p^{w}_{t,i}$ is bounded by the available wind production, i.e., capacity factor, $f_t$ times the wind farm design maximum power $\bar{P}^{w}$:

\begin{equation}
\label{wind_gen_1}
p^{w}_{t,i} \leq f_t \bar{P}^{w}, \quad \forall t \in \mathcal{T},~i \in \mathcal{I}~.
\end{equation}

For simplicity, we only consider one wind production profile from the Ref \cite{rtsgmlc} in this work. 

The operation and maintenance (O\&M) costs $C_{OM}$ for the wind farm and the battery storage are proportional to the maximum wind power $\bar{P}^{w}$ and maximum state of the charge (SoC) $\bar{S}$, respectively:
\begin{equation}
\label{cost_func_1}
C_{OM} = C_{OM}^{w}\bar{P}^{w} + C_{OM}^{b}\bar{S}~.
\end{equation}

Here, $C_{OM}^{w}$ and $C_{OM}^{b}$ represent the operating cost of the wind farm, per unit power (\$/MW), and operation cost of the battery per unit energy (\$/MWh).

We consider a simple linear model for battery storage from Sorourifar et al.\cite{sorourifar2018integrated}, which includes an energy balance to track the SoC of the battery, $S_{t,s}$:
\begin{equation}
\label{batt_soc}
S_{t,i} - S_{t-1,i} = (\eta^{c}p^{c}_{t,i} - p^{d}_{t,i}/\eta^{d})\Delta t, \quad \forall t \in \mathcal{T},~i \in \mathcal{I}.
\end{equation}

The battery charging and discharging efficiency $\eta^{c}$ and $\eta^{d}$ are fixed to 0.95 \cite{armenta2024control, farhad2019introducing, zheng2015study}. The time resolution $\Delta t $ is set to 1 hour to match the PCM's real-time market clearing and wind data resolution. The charging and discharging power are bounded by the battery power capacity $\bar{P}^{b}$.
\begin{equation}
\label{batt_bound}
p^{c}_{t,i} \leq \bar{P}^{b}, \quad p^{d}_{t,i} \leq \bar{P}^{b}, \quad \forall t \in \mathcal{T},~i \in \mathcal{I}.
\end{equation}

We also enforce a periodic boundary condition to account for time beyond the end of the planning horizon
\begin{equation}
\label{batt_init_end}
S_{0, i} = S_{\mathcal{T},i}, \quad i \in \mathcal{I},
\end{equation}
\noindent where $S_{0, i}$ is the initial SoC and $S_{\mathcal{T},i}$ is the final SoC.

To model degradation, we calculate the accumulated energy throughput $E_{t,i}$ through time $t$, the initial energy throughput $E_{0, i}$ is set to $0$:
\begin{equation}
\label{batt_throughput}
E_{t,i} - E_{t-1,i} = \frac{1}{2}(p^{c}_{t,i} + p^{d}_{t,i})\Delta t, \quad \forall t \in \mathcal{T}, \quad i \in \mathcal{I}~.
\end{equation}
The battery SoC $S_{t,i}$ at time $t$ is bounded by the maximum battery SoC $\bar{S}$ minus degradation, which is proportional to the accumulated energy throughput $E_{t,i}$ and the degradation coefficient $\delta$. \cite{sorourifar2018integrated}:
\begin{equation}
\label{batt_degradation}
S_{t,i} \leq \bar{S} - \delta E_{t,i}, \quad \forall t \in \mathcal{T}, \quad i \in \mathcal{I}~.
\end{equation}
We set $\delta$ to 10$^{-4}$, which corresponds to a 50\% loss of full capacity in 5000 cycles.

The maximum battery SoC $\bar{S}$ is the product of the maximum battery power $\bar{P}^{b}$ and the battery capacity of storage $H^{b}$ with the unit of hours.
\begin{equation}
\label{batt_capacity}
\bar{S} = \bar{H} \bar{P}^{b}
\end{equation}

At each time period, the power generated by the wind farm is the sum of power $p^{c}_{t,i}$ charging the battery and power $ p^{s}_{t,i}$ directly selling to the grid:

\begin{equation}
\label{splitter}
p^{w}_{t,i} = p^{c}_{t,i} + p^{s}_{t,i},  \quad \forall t \in \mathcal{T}, \quad i \in \mathcal{I}~.
\end{equation}
The total amount of power provided to the grid by the IES $p_{t,s}$ is the sum of the power discharged from the battery ($p^{d}_{t,s}$) and generated by the wind farm ($p^{s}_{t,s}$):
\begin{equation}
\label{mixer}
p_{t,i} = p^{d}_{t,i} + p^{s}_{t,i}, \quad \forall t \in \mathcal{T}, \quad i \in \mathcal{I}~.
\end{equation}

\subsection{Price-Taker Optimization Model}

The PT co-optimization model determines the optimal storage size and operating decisions to maximize the 30-year NPV:
\begin{subequations}
\begin{align}
\max \quad & \phi \left(\frac{1}{|\mathcal{I}|}\sum_{i \in \mathcal{I}} \sum_{t \in \mathcal{T}} R^{RT}_{t,i} - C_{OM} \right)- C_{capex} \label{pt_model_obj}\\
\mathrm{s.t.} \quad & \text{Eq. } \eqref{wind_gen_1}  \text{  to  }  \text{Eq. } \eqref{mixer}
\label{pt_const1}\\
& R^{RT}_{t,i} = (\hat{\pi}_{t,i}^{RT} + \epsilon)p_{t,i}\Delta t, ~~ \forall t \in \mathcal{T}, ~ i \in \mathcal{I}
\label{pt_const2}\\
& C_{capex} = C^{b}_{capex} \bar{P}^{b}
\label{pt_capex}
\end{align}
\end{subequations}

The NPV discount factor $\phi$ is calculated as follows:
\begin{equation}
\label{NPV}
    \phi = \frac{(1+r)^N - 1}{r(1+r)^N}
\end{equation}

Here, $N = 30$ denotes the life cycle of the IES investment. The discount rate $r$ is set to 0.05.
We consider hourly timesteps over an entire calendar year, i.e., $\mathcal{T} = \{0, 1, 2,..., 8784\}$. We assume that the electricity prices are deterministic, which is modeled with one scenario, i.e., $\mathcal{I} = \{0\}$. Constraint \eqref{pt_const2} calculates the electricity revenue of the IES at time step $t$. $R^{RT}_{t,0}$, denoting the revenue from the real-time market, equals the product of LMP $(\hat{\pi}_{t,0} + \epsilon)$, and power output $p_{t,0}\Delta t$. Here $\epsilon = 10^{-3}$ \$/MWh provides a small incentive to dispatch renewable energy to the grid and eliminates the model degeneracy. Constraint \eqref{pt_capex} calculates the capital cost of the battery storage, $C_{capex}$. All operating and capital costs are from Ref \cite{NRELATB}. 

\subsection{Multiscale Optimization}
We adopt the multiscale optimization (MO) framework from Gao et al.\cite{GAO2022119017} to evaluate the wind-battery IES design candidates. The framework is built based on {\tt Prescient} \cite{prescient}, an open-source PCM package implemented in Python and {\tt Pyomo} \cite{hart2017pyomo}. The {\tt IDAES-PSE} \cite{lee2021idaes} platform provides integration capabilities between {\tt Prescient} and IES bidding and optimization. The Reliability Test System - Grid Modernization Lab Consortium (RTS-GMLC)  network \cite{rtsgmlc, barrows2019ieee} is an open-source dataset developed by the National Renewable Energy Laboratory to expedite power system modeling and grid simulation studies with approximate characteristics of the Southwest United States. 
Using the MO framework, IESs/generators are able to provide time-variant bids to the market. However, in the Prescient PCM simulation, IESs/generators only submit time-invariant bids. In the following subsection, we propose a time-variant bidding method using two-stage stochastic programming. 

\subsection{Time-Variant Bidding}
\label{sb}
We solve the following time-variant bidding problem to compute time-varying bid curves for the wind-battery IES interacting with the market as a dispatchable resource:

\begin{subequations}
\begin{align}
\max ~~ & \frac{1}{|\mathcal{I}|}\sum_{i \in \mathcal{I}}\sum_{t \in \mathcal{T}}\pi^{RT}_{t, i}p^{RT}_{t,i}\Delta t
\label{rt_bidding_obj}\\
\mathrm{s.t.} ~~ & \text{Eq. } \eqref{wind_gen_1} \text{  to  } \text{Eq. }\eqref{mixer} 
\label{rt_bidding_const1}\\
& p^{RT}_{t,i} = p_{t,i}, ~~ \forall t \in \mathcal{T}, \forall i \in \mathcal{I}
\label{rt_bidding_const2}\\
& (p^{RT}_{t,i} - p^{RT}_{t,i'})(\pi^{RT}_{t,i} - \pi^{RT}_{t,i'}) \geq 0, \notag \\
& \qquad \forall i, \forall i^{'} \in \mathcal{I} \backslash i, \forall t \in \mathcal{T}
\label{rt_bidding_const3}
\end{align}
\end{subequations}

The objective function \eqref{rt_bidding_obj} calculates the expected revenue under price uncertainty. Although the IES only participates in the real-time market, we set the planning horizon to four hours, i.e., $\mathcal{T} = \{0, 1, 2, 3, 4\}$. The number of LMP scenarios is set to 10, $\mathcal{I} = \{1, 2, 3,..., 10\}$, which were generated using a  ``backcaster'' and the immediately preceding historical LMPs \cite{elmore2021learning}. Constraint \eqref{rt_bidding_const2} enforces that the real-time offering power is the same as the IES power output to the grid. Constraint \eqref{rt_bidding_const3} enforces the market rule that the offering power is non-decreasing (convex) with respect to marginal cost.

\subsection{Benchmarks and Case Studies}

\begin{table*}[b]
\caption{Comparison of market outcomes from benchmarks and cases (PT: Price-Taker, MO: Multiscale Optimization, TI: Time-Invariant, TV: Time-Variant)}
\label{tb1}
\begin{tabularx}{\textwidth}{@{} l *{10}{C} c @{}}
\toprule
 
                                                  & Benchmark A & Benchmark B  & Benchmark C & Case 1  & Case 2  & Case 3   & Case 4  \\ 
\midrule
IES structure                                     & \multicolumn{3}{c}{Wind Farm Only} &  \multicolumn{4}{c}{Wind Farm + Battery Storage} \\
Battery size ($\bar{H}$, hr)                     & --        & --         & --        & 2       & 10      & 2        & 10 \\
Battery power ratio ($\bar{P}^{b}/\bar{P}^{w}$, MW/MW) & --        & --         & --        & 0.1       & 1.0      & 0.1        & 1.0 \\
Optimization model                                & PT       & PCM         & MO        & PT       & PT      & MO        & MO \\
Bidding strategy                                  & None     & TI         & TV        & None       & None      & TV        & TV \\
Total sold power (GWh)                            & 2089.7     & 1486.5      & 860.6      & 2082.5 & 1990.2 & 912.3  & 1373.8  \\
Total renewable curtailment (GWh)                 & 0        & 603.2       & 1229.1     & 0        & 0      & 1171.0 & 652.3 \\
Battery charging/discharging lost (GWh)           & --        & --         & --        & 7.2      & 99.5   & 6.4  & 63.6 \\
30-yr NPV (M\$)                                   & -251.3    & -251.5   & -181.7   & -259.5   & -2265.9 & -240.0 & -3255.1 \\  
O\&M cost (M\$-yr)                                & 353.9     &  353.9   & 353.9    & 369.7    & 943.5  &  369.7  & 943.5 \\
Capital Cost (M\$)                                & 0         & 0        & 0        & 631.9    & 2358.7 & 631.9    & 2358.7 \\ 
Total electricity revenue (M\$)                   & 19.04       & 19.03        & 23.57       & 24.20   & 100.38  & 25.46   & 36.04 \\ 
Average LMP at the bus (\$/MWh)                   & 24.56       & 24.56        & 24.55       & 24.56   & 24.56   & 25.30    & 25.00  \\ 
Average LMP wind farm/IES received (\$/MWh)       & 9.11        & 12.80        & 27.39       & 11.58   & 50.44   & 29.91    & 26.23  \\ 
Hours that LMP $\geq$ 100 \$/MWh                   & 170         & 170          & 105         & 170     & 170    & 118      & 114    \\ 
Electricity sold when LMP $\geq$ 100 \$/MWh (GWh)  & 13.8       & 13.8         & 8.3        & 21.4   & 148.0  & 8.9    & 14.3   \\ 
\bottomrule
\end{tabularx}
\end{table*}

Three benchmarks and four representative case studies were considered, as summarized in Table \ref{tb1}. This includes comparisons of the original wind farm without any storage using three modeling and bidding strategies: Benchmark A uses PT optimization, Benchmark B uses PCM with TI bids (zero marginal cost), and Benchmark C uses MO with TV bids. Benchmark B generates the reference price signal used for all PT optimization analyses. Cases 1 and 3 consider a wind-battery IES with small battery storage ($\bar{P}^{b}$ = 84.7 MW and $\bar{H}$ = 2 hr), whereas Cases 2 and 4 consider a wind-battery IES with large battery storage ($\bar{P}^{b}$ = 847 MW and $\bar{H}$ = 10 hr). Cases 1 and 2 consider PT optimization, whereas Cases 3 and 4 consider MO with TV bids.


We evenly sampled the design space of battery storage sizes to systematically compare these cases with the benchmarks. The maximum battery power was varied such that the ratio of the maximum battery power ($\bar{P}^{b}$, variable) to wind farm maximum power ($\bar{P}^{w}$, constant) was 0.1 to 1.0 in increments of 0.1. The battery capacity $\bar{H}$ was varied from 2 to 10 hours with increments of 2 hours.

{\tt Prescient} (version 2.2.2) and {\tt Gurobi} (version 9.5.1) \cite{gurobi} with a 1\% relative mixed-integer programming gap were used to solve the unit commitment problems. {\tt Ipopt} (version 3.14.2) \cite{wachter2006implementation} with HSL linear algebra routines \cite{HSL} was used to solve the time-variant bidding and PT optimization problems.

\IEEEpubidadjcol

\section{Results and Discussion}

Table \ref{tb1} summarizes the results. Additional results and discussion are organized around three key findings.

\label{sec:3}
\subsection{Price-Taker is Overly Optimistic}

\begin{figure}[h]
\centering
\includegraphics[width=3.5in]{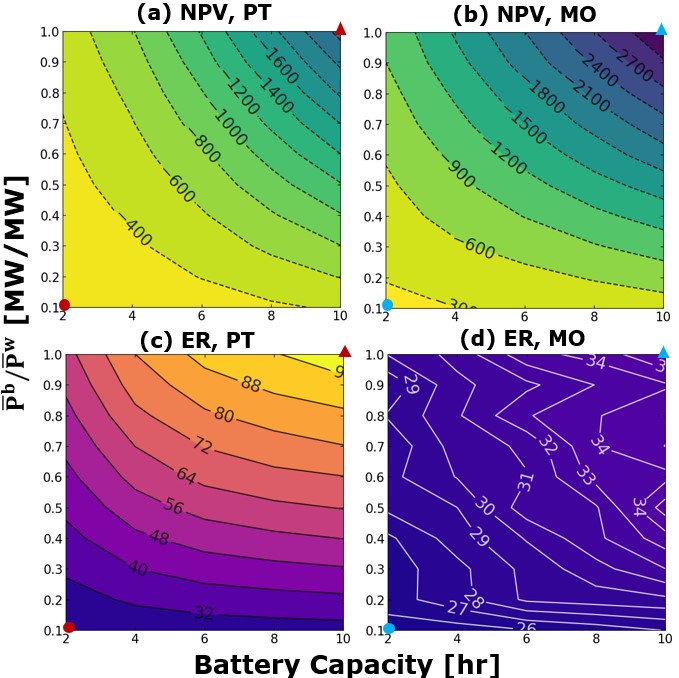}
\caption{Comparison of net present value, NPV (M\$, \textbf{(a) and (b)}) and annual electricity revenue, ER (M\$, \textbf{(c) and (d)}) results from PT (\textbf{(a) and (c)}) and MO (\textbf{(b) and (d)}). The symbols correspond to scenarios in Table \ref{tb1}: (red dots) Case 1, (red triangles) Case 2, (blue dots) Case 3, and (blue triangles) Case 4.}
\label{fig:npv_elecrev_comparison}
\end{figure}

We note that all reported NPVs are negative. However, the PT optimization overestimates the electricity revenue from the market compared to the MO, leading to inaccurate NPV predictions and conclusions. Fig. \ref{fig:npv_elecrev_comparison} (a) and (b) compare the 30-year NPV of the IES retrofit. Despite the similar trends, the NPV estimated by the PT is higher than that of MO. From Fig. \ref{fig:npv_elecrev_comparison} (c) and (d), the PT result shows IES electricity revenue increased by 314.9\% from adding the smallest battery (Case 1, red dot) to the largest battery (Case 2, red triangle), while in the MO, it only increased 41.1\%. From the PT optimization, the maximum electricity revenue of the wind-battery IES is up to 100.38 M\$ (Case 2, red triangle). However, in the MO, although the maximum electricity revenue is still obtained in the maximum battery storage design (Case 4, blue triangle), the value is only 36.04 M\$, which is only 35.9\% of the result from PT optimization. The PT is overly optimistic because it assumes all electricity can be sold at given price signals, thus overestimating the electricity revenue and NPV.

\begin{figure}[h]
\centering
\includegraphics[width=3.5in]{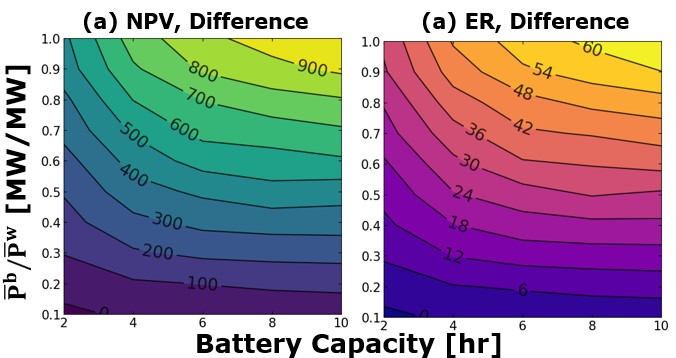}
\caption{Difference (PT - MO) of net present value, NPV (M\$, \textbf{(a)}) and annual electricity revenue, ER (M\$, \textbf{(b)}). The difference is more sensitive to the $\bar{P}^{b}/\bar{P}^{w}.$} 
\label{fig:npv_elecrev_diff}
\end{figure}

Fig. \ref{fig:npv_elecrev_comparison} reveals that electricity revenue is more sensitive to the maximum battery power $\bar{P}^{b}$. The RTS-GMLC network has excess generation and transmission capacity for the demand profiles \cite{barrows2019ieee}. Thus, periods of high electricity prices are short, and electricity prices are low; therefore, larger $\bar{P}^{b}$ enables the IES to better exploit short periods of high prices. However, the PT optimization overestimates revenue and NPV more with larger $\bar{P}^{b}$. Fig. \ref{fig:npv_elecrev_diff} visualizes the differences between PT and MO results, which confirms these trends. Because RTS-GMLC is overbuilt and prices are low, the MO selects the smallest allowed battery, and this design has a negative NPV value, which means the IES retrofit is not economical in RTS-GMLC as a dispatchable resource.

The results reported in Table \ref{tb1} show PT overestimates the renewable integration capability of the wind-battery IES compared to MO. Specifically, the PT approach (Benchmark A, Case 1, and Case 2) predicts significant energy losses from the battery round trip efficiency and no renewable curtailment. 
In Benchmark B (PCM and TI bids), however, the IES bids the power at zero marginal cost, but 28.9\% wind energy is curtailed due to grid constraints. Thus PT overestimates the renewable integration. With MO (Cases 3 and 4), the renewable curtailment is 1171.0 GWh and 652.3 GWh higher than with PT (Cases 1 and 2). This result makes sense because the PT approach neglects grid constraints and assumes that all the electricity from the IES can be sold to the grid, whereas MO explicitly considers grid constraints.

In summary, these results caution that the PT approach overestimates the IES retrofit's economic value and renewable integration capability, leading to overly optimistic techno-economic analyses.   

\subsection{Strategic Bidding Improves IES Revenues}
Next, we explore how time-variant bids impact the IES operation and market clearing. We start by considering only the wind farm. In Benchmark A (PT), all electricity generated from the wind farm (2089.7 GWh) is sold to the market due to the PT assumption. In contrast, with PCM simulation, which considers more realistic grid constraints (Benchmark B), only 71.1\% of the available wind energy is sold to the grid. 
However, Benchmarks A and B have almost identical annual electricity revenues. This is because in Benchmark A the additional electricity is sold to the grid at times when the price is nearly zero.

\begin{figure}[h]
\centering
\includegraphics[width=3.5in]{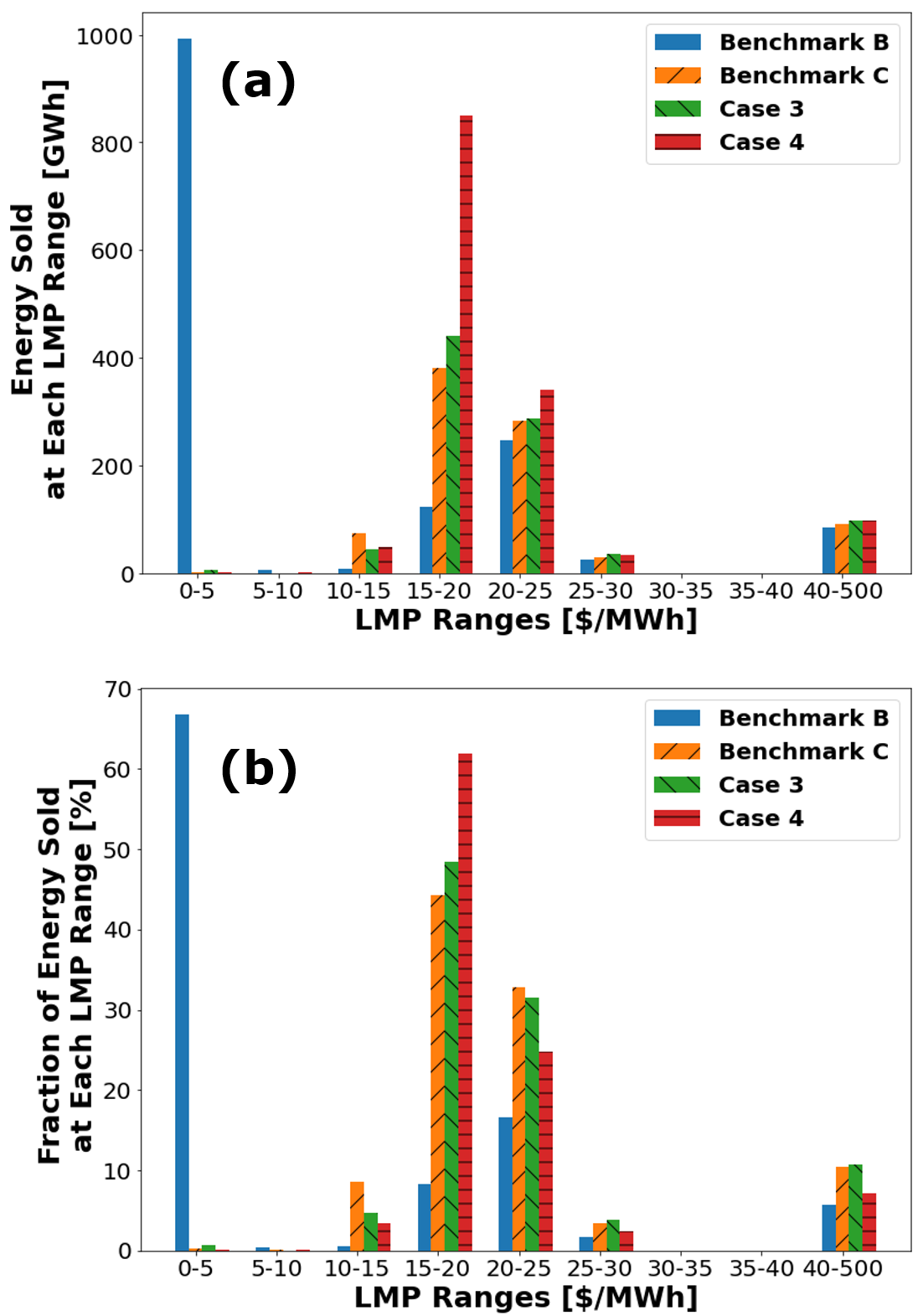}
\caption{Comparison of energy sold at different LMP ranges (a) and its fraction (b). Using the time-variant bidding, the wind generator has much less power dispatched under 5 \$/MWh. Also, the battery power is dispatched mostly between 15 \$/MWh - 25 \$/MWh.}
\label{fig:power_zone}
\end{figure}

With time-varying (TV) bids for the wind farm (Benchmark C), the electricity revenue increases by 23.9\% while only selling 41.2\% and 57.9\% of the electricity relative to Benchmarks A and B, respectively. This confirms that TV bids benefit the resource. Moreover, the average LMP the wind farm received was 27.39 \$/MWh, more than twice that in Benchmark B, but the average LMP at the bus was nearly constant across Benchmarks A to C. This suggests that the TV bids do not manipulate the market prices because the wind farm is not the marginal unit within the grid. Fig. \ref{fig:power_zone} shows the energy dispatched as a function of LMP. In Benchmark B (PCM, TI bids), 66.8\% of the electricity is sold at prices between 0 and 5 \$/MWh. In Benchmark C (MO, TV bids), 0.3\% of the energy is sold at prices under 5 \$/MWh, and 77.2\% of the energy is sold at prices between 15 and 25 \$/MWh. The result is consistent with the ``merit-order effect'' discussed in Refs \cite{gomes2023hybrid, antweiler2021long}.

\begin{figure}[h]
\centering
\includegraphics[width=3.5in]{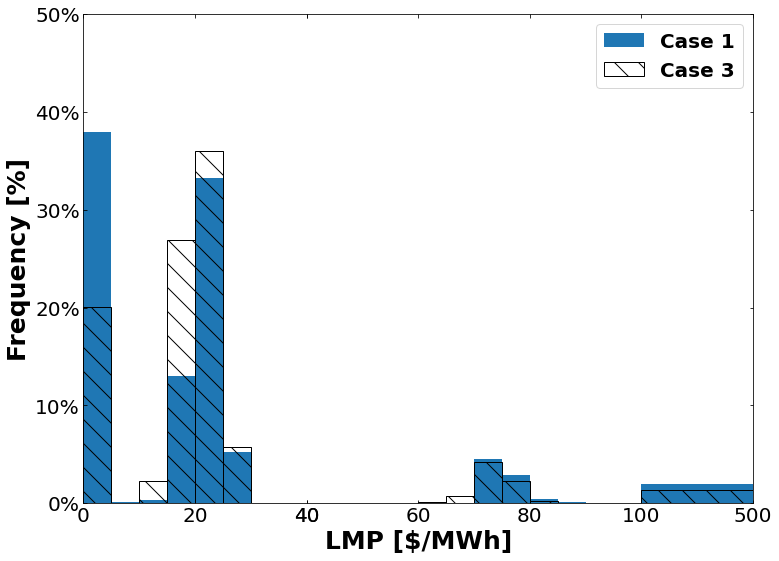}
\caption{The histogram of annual real-time LMP distribution for scenarios considering bidding. From Case 1 to Case 3, hours of LMP lower than 5 \$/MWh decreased from 37.9\% to 20.1\%. Hours of LMP between 15 \$/MWh and 25 \$/MWh are increased from 46.24\% to 62.92\%. } 
\label{fig:lmp_hist}
\end{figure}

Next, we compare the wind-battery IES using TV bidding (Cases 3 and 4) with the price-taker approach (Cases 1 and 2). In Case 1 (PT, small battery) and Case 2 (PT, large battery), the IES sold 21.4 GWh and 148.0 GWh energy, respectively, when the LMP was $\geq$ 100 \$/MWh, which is referred to as the ``high LMP hour." However, in Case 3 (MO, small battery) and Case 4 (MO, large battery), the IES only sold 41.6\% and 9.7\% power at the ``high LMP hour" compared to Cases 1 and 2. This indicates that while strategic bidding can improve the IES revenue, it cannot exploit all the high-price opportunities due to the grid constraints. This is further discussed in the next subsection. Moreover, the number of ``high LMP hours" in Cases 3 and 4 is 69.4\% and 67.1\% of that in Cases 1 and 2, respectively, showing the IES with strategic bidding impacts market prices by providing more flexibility to the grid, thus reducing the number of "high LMP hours." This nuanced interaction cannot be predicted by the price-taker model. Fig. \ref{fig:lmp_hist} shows this finding by comparing the annual LMP distribution at the ``Caesar" bus in Cases 1 and 3. The IES-market interaction reduces the number of periods that the LMP goes to zero and increases periods of LMP between 15 to 25 \$/MWh.    

Finally, we compare Benchmark C (MO, TV bids) with Cases 3 and 4 to demonstrate the benefit of battery storage. A small battery (Case 3) and a large battery (Case 4) reduced the renewable curtailment by 58.1 GWh and 576.8 GWh, respectively, and increased electricity revenue by 8.0\% (Case 3) and 52.9\% (Case 4) compared with Benchmark C. However, these revenues are not sufficient to offset the cost of storage and achieve positive NPVs. This motivates future work to explore additional revenue opportunities from ancillary service markets, which prior work showed (using PT) can significantly boost revenues and result in positive NPVs for battery energy storage\cite{sorourifar2018integrated}.       

Thus in summary, with TV bids, the wind-battery IES becomes the ``price-maker'' and significantly influences market prices. PT optimization seeks to dispatch as much electricity as possible, maximizing its revenue based on (perfect) knowledge of all price signals. 
In contrast, MO is a more nuanced approach that considers IES-market interactions including that the market will rarely accept all of the power production capacity offered by the IES. Moreover, MO considers uncertainty in price signals when constructing TV bids.

\subsection{Price-Taker Neglects the Price Uncertainty.}
\label{sec:3_3}
\begin{figure}[h]
\centering
\includegraphics[width=3.5in]{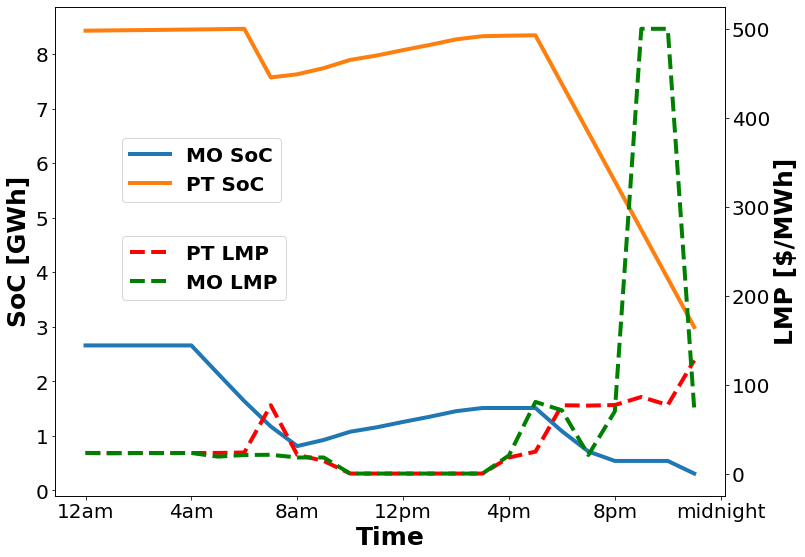}
\caption{01-19-2020 simulation results in Case 4 exhibit the difference between the battery SoC and LMPs in the PT and MO.}
\label{fig:soclmp}
\end{figure}

PT and MO also predict different optimal IES operating decisions, as illustrated in the example in Fig. \ref{fig:soclmp}. During the simulation date 01-19-2020, high LMP hours in both models start at 5 pm. In PT optimization, due to the perfect information and simplistic grid abstraction, the IES charges the battery before 5 pm, and the SoC reaches 8344 MWh. The battery keeps discharging at the maximum rate from hour 5 pm to hour 11 pm. With MO, 9 pm and 10 pm have two price spikes (500 \$/MWh), which are not forecasted correctly. The IES at 9 pm and 10 pm bids all the available wind power but does not bid for any battery power. Also, the SoC level is relatively lower (1512 MWh) at 5 pm. As a result, in PT, the IES has 595k\$ revenue while the MO only has 352k\$ on this day. The PT model (incorrectly) assumes perfect information about the LMP signal. We note that other prior studies quantified the impact of perfect information by considering LMP forecasts in PT.\cite{elmore2021learning} Nevertheless, MO is more representative of real-world market operations, and the TV bidding considers LMP uncertainty.

\section{Conclusion}
\label{sec:4}
This work compares techno-economic optimization strategies for retrofitting wind generators with a battery storage system to create a dispatchable IES. We show PT optimization overestimates the IES annual electricity revenue up to 178\%, resulting in a maximum NPV overestimation of 50\%. This is because the PT approximation assumes LMPs are known \emph{a priori} and not significantly impacted by the IES. Using a rigorous MO framework, we show this approximation is inaccurate. PT approximation also overestimates the capability of renewable integration by assuming the grid can accept all the electricity from the IES. Moreover, we find that TV bids increase wind farm electricity revenue by 24\% and the selling price of electricity by 114\%. TV bids make the wind-battery IES the ``price-maker,'' i.e., marginal generator, and thus its actions dramatically impact market prices. Thus, while the PT assumption is well-suited for preliminary analysis, we emphasize the need to go beyond PT when optimizing the economics of IESs.

Several improvements are recommended for future work. We anticipate more accurate LMP forecasting methods\cite{elmore2021learning} will further improve the performance of the MO framework. Sophisticated bidding strategies beyond those in this paper should be benchmarked using the MO framework instead of PT to ensure IES-market interactions are considered. The retirement and replacement of battery storage caused by the degradation should also be considered in the decision-making process. Finally, scale-bridging surrogate models can be used to incorporate more complex market interactions revealed by PCM simulations to co-optimize the design, operation, and market participation of renewable-storage IESs. 

\section*{Credit Statement}
\textbf{Xinhe Chen:} Conceptualization, Methodology, Software, Formal analysis, Investigation, Data curation, Writing– original draft, Writing– review \& editing, Visualization. 
\textbf{Xian Gao:} Conceptualization, Methodology, Software, Formal analysis, Investigation, Writing– original draft. \textbf{Darice Guittet:} Conceptualization, Methodology, Software, Formal analysis, Investigation, Writing– original draft. \textbf{Radhakrishna Tumbalam-Gooty:} Conceptualization, Methodology, Software. 
\textbf{Bernard Knueven:} Conceptualization, Methodology, Software. 
\textbf{John D. Siirola:} Project administration, Conceptualization, Methodology. 
\textbf{David C. Miller:} Project administration, Funding acquisition, Conceptualization.
\textbf{Alexander W. Dowling:} Project administration, Conceptualization, Methodology, Software, Formal analysis, Investigation, Data curation, Writing– review \& editing, Visualization.

\section*{Acknowledgments}
This work was conducted as part of the Design Integration and Synthesis Platform to Advance Tightly Coupled Hybrid Energy Systems (DISPATCHES) project through the Grid Modernization Lab Consortium with funding from the U.S. Department of Energy's Office of Fossil Energy and Carbon Management, Office of Nuclear Energy, and Hydrogen and Fuel Cell Technologies Office. Additional work was conducted as part of the Institute for the Design of Advanced Energy Systems (IDAES) with support from the U.S. Department of Energy's Office of Fossil Energy and Carbon Management (FECM) through the Simulation-based Engineering Research Program. 

We thank Jaffer Ghouse and Wesley Jones for technical input and help with open-source software release management, and thank Keith Beattie and Ludovico Bianchi for the help with open-source software release management and coordination.

\textbf{Disclaimer}
This project was funded by the Department of Energy, National Energy Technology Laboratory an agency of the United States Government, through a support contract. Neither the United States Government nor any agency thereof, nor any of its employees, nor the support contractor, nor any of their employees, makes any warranty, expressor implied, or assumes any legal liability or responsibility for the accuracy, completeness, or usefulness of any information, apparatus, product, or process disclosed, or represents that its use would not infringe privately owned rights. Reference herein to any specific commercial product, process, or service by trade name, trademark, manufacturer, or otherwise does not necessarily constitute or imply its endorsement, recommendation, or favoring by the United States Government or any agency thereof. The views and opinions of authors expressed herein do not necessarily state or reflect those of the United States Government or any agency thereof. 

This work was authored [in part] by the National Renewable Energy Laboratory, operated by Alliance for Sustainable Energy, LLC, for the U.S. Department of Energy (DOE) under Contract No. DE-AC36-08GO28308. The views expressed in the article do not necessarily represent the views of the DOE or the U.S. Government. The U.S. Government retains and the publisher, by accepting the article for publication, acknowledges that the U.S. Government retains a nonexclusive, paid-up, irrevocable, worldwide license to publish or reproduce the published form of this work, or allow others to do so, for U.S. Government purposes.

Sandia National Laboratories is a multi-mission laboratory managed and operated by National Technology \& Engineering Solutions of Sandia, LLC (NTESS), a wholly owned subsidiary of Honeywell International Inc., for the U.S. Department of Energy’s National Nuclear Security Administration (DOE/NNSA) under contract DE-NA0003525. This written work is authored by an employee of NTESS. The employee, not NTESS, owns the right, title and interest in and to the written work and is responsible for its contents. Any subjective views or opinions that might be expressed in the written work do not necessarily represent the views of the U.S. Government. The publisher acknowledges that the U.S. Government retains a non-exclusive, paid-up, irrevocable, world-wide license to publish or reproduce the published form of this written work or allow others to do so, for U.S. Government purposes. The DOE will provide public access to results of federally sponsored research in accordance with the DOE Public Access Plan.




\begin{thebibliography}{10}
\providecommand{\url}[1]{#1}
\csname url@samestyle\endcsname
\providecommand{\newblock}{\relax}
\providecommand{\bibinfo}[2]{#2}
\providecommand{\BIBentrySTDinterwordspacing}{\spaceskip=0pt\relax}
\providecommand{\BIBentryALTinterwordstretchfactor}{4}
\providecommand{\BIBentryALTinterwordspacing}{\spaceskip=\fontdimen2\font plus
\BIBentryALTinterwordstretchfactor\fontdimen3\font minus \fontdimen4\font\relax}
\providecommand{\BIBforeignlanguage}[2]{{%
\expandafter\ifx\csname l@#1\endcsname\relax
\typeout{** WARNING: IEEEtran.bst: No hyphenation pattern has been}%
\typeout{** loaded for the language `#1'. Using the pattern for}%
\typeout{** the default language instead.}%
\else
\language=\csname l@#1\endcsname
\fi
#2}}
\providecommand{\BIBdecl}{\relax}
\BIBdecl

\bibitem{worldrenewable}
H.~Ritchie, M.~Roser, and P.~Rosado, ``Renewable energy,'' \emph{Our World in Data}, 2020, https://ourworldindata.org/renewable-energy.

\bibitem{unrenewable}
\BIBentryALTinterwordspacing
``{United Nations} climate action,'' accessed: 2024-06-13. [Online]. Available: \url{https://www.un.org/en/climatechange/raising-ambition/renewable-energy}
\BIBentrySTDinterwordspacing

\bibitem{usdoe}
\BIBentryALTinterwordspacing
``{U.S.}~installed and potential wind power capacity and generation,'' accessed: 2024-04-22. [Online]. Available: \url{https://windexchange.energy.gov/maps-data/321}
\BIBentrySTDinterwordspacing

\bibitem{vargas2019wind}
S.~A. Vargas, G.~R.~T. Esteves, P.~M. Ma{\c{c}}aira, B.~Q. Bastos, F.~L.~C. Oliveira, and R.~C. Souza, ``Wind power generation: A review and a research agenda,'' \emph{Journal of Cleaner Production}, vol. 218, pp. 850--870, 2019.

\bibitem{ertugrul2016battery}
N.~Ertugrul, ``Battery storage technologies, applications and trend in renewable energy,'' in \emph{2016 IEEE International Conference on Sustainable Energy Technologies (ICSET)}.\hskip 1em plus 0.5em minus 0.4em\relax IEEE, 2016, pp. 420--425.

\bibitem{arent2021multi}
D.~J. Arent, S.~M. Bragg-Sitton, D.~C. Miller, T.~J. Tarka, J.~A. Engel-Cox, R.~D. Boardman, P.~C. Balash, M.~F. Ruth, J.~Cox, and D.~J. Garfield, ``Multi-input, multi-output hybrid energy systems,'' \emph{Joule}, vol.~5, no.~1, pp. 47--58, 2021.

\bibitem{kang2011optimal}
C.~A. Kang, A.~R. Brandt, and L.~J. Durlofsky, ``Optimal operation of an integrated energy system including fossil fuel power generation, co2 capture and wind,'' \emph{Energy}, vol.~36, no.~12, pp. 6806--6820, 2011.

\bibitem{wang2019operation}
Y.~Wang, Y.~Wang, Y.~Huang, J.~Yang, Y.~Ma, H.~Yu, M.~Zeng, F.~Zhang, and Y.~Zhang, ``Operation optimization of regional integrated energy system based on the modeling of electricity-thermal-natural gas network,'' \emph{Applied Energy}, vol. 251, p. 113410, 2019.

\bibitem{chen2020optimal}
X.~Chen, C.~Wang, Q.~Wu, X.~Dong, M.~Yang, S.~He, and J.~Liang, ``Optimal operation of integrated energy system considering dynamic heat-gas characteristics and uncertain wind power,'' \emph{Energy}, vol. 198, p. 117270, 2020.

\bibitem{bragg2020reimagining}
S.~M. Bragg-Sitton, R.~Boardman, C.~Rabiti, and J.~O'Brien, ``Reimagining future energy systems: Overview of the us program to maximize energy utilization via integrated nuclear-renewable energy systems,'' \emph{International Journal of Energy Research}, vol.~44, no.~10, pp. 8156--8169, 2020.

\bibitem{li2017optimal}
G.~Li, R.~Zhang, T.~Jiang, H.~Chen, L.~Bai, H.~Cui, and X.~Li, ``Optimal dispatch strategy for integrated energy systems with cchp and wind power,'' \emph{Applied energy}, vol. 192, pp. 408--419, 2017.

\bibitem{li2021two}
P.~Li, Z.~Wang, J.~Wang, W.~Yang, T.~Guo, and Y.~Yin, ``Two-stage optimal operation of integrated energy system considering multiple uncertainties and integrated demand response,'' \emph{Energy}, vol. 225, p. 120256, 2021.

\bibitem{li2020improving}
Y.~Li, C.~Wang, G.~Li, J.~Wang, D.~Zhao, and C.~Chen, ``Improving operational flexibility of integrated energy system with uncertain renewable generations considering thermal inertia of buildings,'' \emph{Energy Conversion and Management}, vol. 207, p. 112526, 2020.

\bibitem{jiang2020optimal}
P.~Jiang, J.~Dong, and H.~Huang, ``Optimal integrated demand response scheduling in regional integrated energy system with concentrating solar power,'' \emph{Applied Thermal Engineering}, vol. 166, p. 114754, 2020.

\bibitem{ruiming2019multi}
F.~Ruiming, ``Multi-objective optimized operation of integrated energy system with hydrogen storage,'' \emph{International Journal of Hydrogen Energy}, vol.~44, no.~56, pp. 29\,409--29\,417, 2019.

\bibitem{fabianek2023techno}
P.~Fabianek and R.~Madlener, ``Techno-economic analysis and optimal sizing of hybrid pv-wind systems for hydrogen production by pem electrolysis in {California} and {Northern Germany},'' \emph{International Journal of Hydrogen Energy}, 2023.

\bibitem{sioshansi2021energy}
R.~Sioshansi, P.~Denholm, J.~Arteaga, S.~Awara, S.~Bhattacharjee, A.~Botterud, W.~Cole, A.~Cortes, A.~De~Queiroz, J.~DeCarolis \emph{et~al.}, ``Energy-storage modeling: State-of-the-art and future research directions,'' \emph{IEEE transactions on power systems}, vol.~37, no.~2, pp. 860--875, 2021.

\bibitem{koohi2020review}
S.~Koohi-Fayegh and M.~A. Rosen, ``A review of energy storage types, applications and recent developments,'' \emph{Journal of Energy Storage}, vol.~27, p. 101047, 2020.

\bibitem{yang2018battery}
Y.~Yang, S.~Bremner, C.~Menictas, and M.~Kay, ``Battery energy storage system size determination in renewable energy systems: A review,'' \emph{Renewable and Sustainable Energy Reviews}, vol.~91, pp. 109--125, 2018.

\bibitem{denholm2023moving}
P.~Denholm, W.~Cole, and N.~Blair, ``Moving beyond 4-hour li-ion batteries: Challenges and opportunities for long (er)-duration energy storage,'' National Renewable Energy Laboratory (NREL), Golden, CO, Tech. Rep. NREL/TP-6A40-85878, 2023.

\bibitem{dowling2017multi}
A.~W. Dowling, R.~Kumar, and V.~M. Zavala, ``A multi-scale optimization framework for electricity market participation,'' \emph{Applied Energy}, vol. 190, pp. 147--164, 2017.

\bibitem{bechlenberg2024renewable}
A.~Bechlenberg, E.~A. Luning, M.~B. Salt{\i}k, N.~B. Szirbik, B.~Jayawardhana, and A.~I. Vakis, ``Renewable energy system sizing with power generation and storage functions accounting for its optimized activity on multiple electricity markets,'' \emph{Applied Energy}, vol. 360, p. 122742, 2024.

\bibitem{Dicorato2012planning}
M.~Dicorato, G.~Forte, M.~Pisani, and M.~Trovato, ``Planning and operating combined wind-storage system in electricity market,'' \emph{IEEE Transactions on Sustainable Energy}, vol.~3, no.~2, pp. 209--217, 2012.

\bibitem{Moghaddam2019bptimal}
I.~N. Moghaddam, B.~Chowdhury, and M.~Doostan, ``Optimal sizing and operation of battery energy storage systems connected to wind farms participating in electricity markets,'' \emph{IEEE Transactions on Sustainable Energy}, vol.~10, no.~3, pp. 1184--1193, 2019.

\bibitem{hou2018cooperation}
P.~Hou, G.~Yang, P.~Enevoldsen, and A.~H. Nielsen, ``Cooperation of offshore wind farm with battery storage in multiple electricity markets,'' in \emph{2018 53rd international universities power engineering conference (UPEC)}.\hskip 1em plus 0.5em minus 0.4em\relax IEEE, 2018, pp. 1--6.

\bibitem{he2016cooperation}
G.~He, Q.~Chen, C.~Kang, Q.~Xia, and K.~Poolla, ``Cooperation of wind power and battery storage to provide frequency regulation in power markets,'' \emph{IEEE Transactions on Power Systems}, vol.~32, no.~5, pp. 3559--3568, 2016.

\bibitem{nguyen2012new}
M.~Y. Nguyen, D.~H. Nguyen, and Y.~T. Yoon, ``A new battery energy storage charging/discharging scheme for wind power producers in real-time markets,'' \emph{Energies}, vol.~5, no.~12, pp. 5439--5452, 2012.

\bibitem{GAO2022119017}
X.~Gao, B.~Knueven, J.~D. Siirola, D.~C. Miller, and A.~W. Dowling, ``Multiscale simulation of integrated energy system and electricity market interactions,'' \emph{Applied Energy}, vol. 316, p. 119017, 2022.

\bibitem{sorourifar2018integrated}
F.~Sorourifar, V.~M. Zavala, and A.~W. Dowling, ``Integrated multiscale design, market participation, and replacement strategies for battery energy storage systems,'' \emph{IEEE Transactions on Sustainable Energy}, vol.~11, no.~1, pp. 84--92, 2018.

\bibitem{emmanuel2022market}
M.~I. Emmanuel and P.~Denholm, ``A market feedback framework for improved estimates of the arbitrage value of energy storage using price-taker models,'' \emph{Applied Energy}, vol. 310, p. 118250, 2022.

\bibitem{frew2023analysis}
B.~Frew, D.~Levie, J.~Richards, J.~Desai, and M.~Ruth, ``Analysis of multi-output hybrid energy systems interacting with the grid: Application of improved price-taker and price-maker approaches to nuclear-hydrogen systems,'' \emph{Applied Energy}, vol. 329, p. 120184, 2023.

\bibitem{teixeira2012market}
F.~Teixeira, J.~de~Sousa, and S.~Faias, ``How market power affects the behavior of a pumped storage hydro unit in the day-ahead electricity market?'' in \emph{2012 9th International Conference on the European Energy Market}.\hskip 1em plus 0.5em minus 0.4em\relax IEEE, 2012, pp. 1--6.

\bibitem{batlle2005strategic}
C.~Batlle and J.~Barqu{\'\i}n, ``A strategic production costing model for electricity market price analysis,'' \emph{IEEE Transactions on Power Systems}, vol.~20, no.~1, pp. 67--74, 2005.

\bibitem{jalving2023beyond}
J.~Jalving, J.~Ghouse, N.~Cortes, X.~Gao, B.~Knueven, D.~Agi, S.~Martin, X.~Chen, D.~Guittet, R.~Tumbalam-Gooty \emph{et~al.}, ``Beyond price taker: Conceptual design and optimization of integrated energy systems using machine learning market surrogates,'' \emph{Applied Energy}, vol. 351, p. 121767, 2023.

\bibitem{cocchi2018machine}
G.~Cocchi, L.~Galli, G.~Galvan, M.~Sciandrone, M.~Cant{\`u}, and G.~Tomaselli, ``Machine learning methods for short-term bid forecasting in the renewable energy market: A case study in italy,'' \emph{Wind Energy}, vol.~21, no.~5, pp. 357--371, 2018.

\bibitem{zheng2023energy}
N.~Zheng, X.~Qin, D.~Wu, G.~Murtaugh, and B.~Xu, ``Energy storage state-of-charge market model,'' \emph{IEEE Transactions on Energy Markets, Policy and Regulation}, vol.~1, no.~1, pp. 11--22, 2023.

\bibitem{kian2005bidding}
A.~R. Kian, J.~B. Cruz, and R.~J. Thomas, ``Bidding strategies in oligopolistic dynamic electricity double-sided auctions,'' \emph{IEEE Transactions on Power Systems}, vol.~20, no.~1, pp. 50--58, 2005.

\bibitem{rtsgmlc}
\BIBentryALTinterwordspacing
``Reliability test system grid modernization lab consortium,'' accessed: 2024-05-21. [Online]. Available: \url{https://github.com/GridMod/RTS-GMLC}
\BIBentrySTDinterwordspacing

\bibitem{armenta2024control}
C.~Armenta-Deu and T.~Coulaud, ``Control unit for battery charge management in electric vehicles (evs),'' \emph{Future Transportation}, vol.~4, no.~2, pp. 429--449, 2024.

\bibitem{farhad2019introducing}
S.~Farhad and A.~Nazari, ``Introducing the energy efficiency map of lithium-ion batteries,'' \emph{International Journal of Energy Research}, vol.~43, no.~2, pp. 931--944, 2019.

\bibitem{zheng2015study}
Y.~Zheng, M.~Ouyang, L.~Lu, J.~Li, Z.~Zhang, and X.~Li, ``Study on the correlation between state of charge and coulombic efficiency for commercial lithium ion batteries,'' \emph{Journal of power Sources}, vol. 289, pp. 81--90, 2015.

\bibitem{NRELATB}
\BIBentryALTinterwordspacing
``2023 electricity {ATB} technologies,'' accessed: 2023-06-30. [Online]. Available: \url{https://atb.nrel.gov/electricity/2023/technologies}
\BIBentrySTDinterwordspacing

\bibitem{prescient}
\BIBentryALTinterwordspacing
``Prescient production cost modeling platform,'' accessed: 2024-05-21. [Online]. Available: \url{https://github.com/grid-parity-exchange/Prescient}
\BIBentrySTDinterwordspacing

\bibitem{hart2017pyomo}
W.~E. Hart, C.~D. Laird, J.-P. Watson, D.~L. Woodruff, G.~A. Hackebeil, B.~L. Nicholson, J.~D. Siirola \emph{et~al.}, \emph{Pyomo-optimization modeling in python}.\hskip 1em plus 0.5em minus 0.4em\relax Springer, 2017, vol.~67.

\bibitem{lee2021idaes}
A.~Lee, J.~H. Ghouse, J.~C. Eslick, C.~D. Laird, J.~D. Siirola, M.~A. Zamarripa, D.~Gunter, J.~H. Shinn, A.~W. Dowling, D.~Bhattacharyya \emph{et~al.}, ``The idaes process modeling framework and model library—flexibility for process simulation and optimization,'' \emph{Journal of advanced manufacturing and processing}, vol.~3, no.~3, p. e10095, 2021.

\bibitem{barrows2019ieee}
C.~Barrows, A.~Bloom, A.~Ehlen, J.~Ik{\"a}heimo, J.~Jorgenson, D.~Krishnamurthy, J.~Lau, B.~McBennett, M.~O’Connell, E.~Preston \emph{et~al.}, ``The ieee reliability test system: A proposed 2019 update,'' \emph{IEEE Transactions on Power Systems}, vol.~35, no.~1, pp. 119--127, 2019.

\bibitem{elmore2021learning}
C.~T. Elmore and A.~W. Dowling, ``Learning spatiotemporal dynamics in wholesale energy markets with dynamic mode decomposition,'' \emph{Energy}, vol. 232, p. 121013, 2021.

\bibitem{gurobi}
\BIBentryALTinterwordspacing
``Gurobi optimization,'' accessed: 2024-06-17. [Online]. Available: \url{https://gurobi.com}
\BIBentrySTDinterwordspacing

\bibitem{wachter2006implementation}
A.~W{\"a}chter and L.~T. Biegler, ``On the implementation of an interior-point filter line-search algorithm for large-scale nonlinear programming,'' \emph{Mathematical programming}, vol. 106, pp. 25--57, 2006.

\bibitem{HSL}
\BIBentryALTinterwordspacing
``Hsl, a collection of fortran codes for large-scale scientific computation.'' accessed: 2024-07-02. [Online]. Available: \url{http://hsl.rl.ac.uk/}
\BIBentrySTDinterwordspacing

\bibitem{gomes2023hybrid}
J.~G. Gomes, J.~Jiang, C.~T. Chong, J.~Telhada, X.~Zhang, S.~Sammarchi, S.~Wang, Y.~Lin, and J.~Li, ``Hybrid solar pv-wind-battery system bidding optimisation: A case study for the iberian and italian liberalised electricity markets,'' \emph{Energy}, vol. 263, p. 126043, 2023.

\bibitem{antweiler2021long}
W.~Antweiler and F.~Muesgens, ``On the long-term merit order effect of renewable energies,'' \emph{Energy Economics}, vol.~99, p. 105275, 2021.

\end{thebibliography}


\begin{IEEEbiography}[{\includegraphics[width=1in,height=1.25in,clip,keepaspectratio]{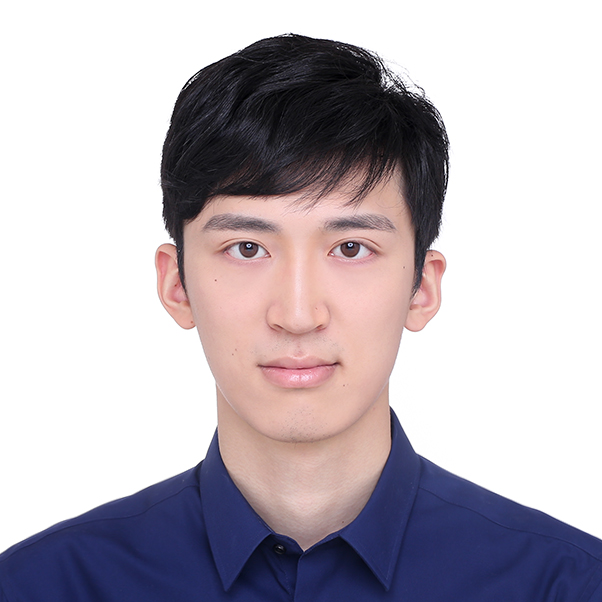}}]{Xinhe Chen}
received a B.S.E degree in Resource Recycling Science and Engineering from East China University of Science and Technology, Shanghai, China and an M.S. degree in Chemical Engineering from Carnegie Mellon University. He is currently a Ph.D student in Chemical Engineering at the University of Notre Dame, Notre Dame, IN. 
\end{IEEEbiography}
\begin{IEEEbiography}[{\includegraphics[width=1in,height=1.25in,clip,keepaspectratio]{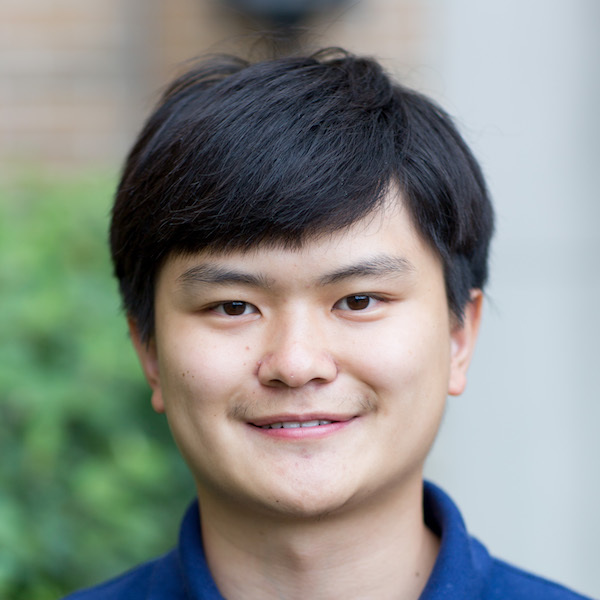}}]{Xian Gao}
received a B.S.E degree in Chemical Engineering from Hunan University, Hunan, China, and a Ph.D degree in Chemical Engineering from the University of Notre Dame, Notre Dame, IN. He is now an operations research scientist at Grubhub Holdings Inc., in Chicago, IL.  
\end{IEEEbiography}

\begin{IEEEbiography}[{\includegraphics[width=1in,height=1.25in,clip,keepaspectratio]{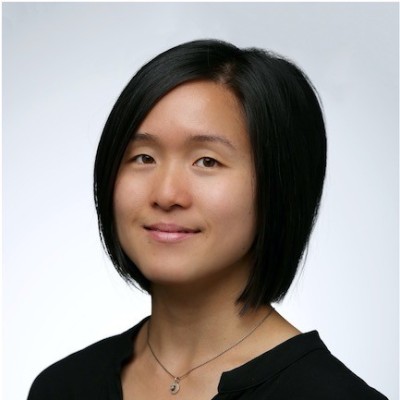}}]{Darice Guittet}
received a B.S. degree in Applied Mathematics from the University of Colorado Boulder. She is currently a researcher in the Strategic Energy Analysis Center's Distributed Systems and Storage Group at the National Renewable Energy Laboratory in Golden, CO. 
\end{IEEEbiography}

\begin{IEEEbiography}[{\includegraphics[width=1in,height=1.25in,clip,keepaspectratio]{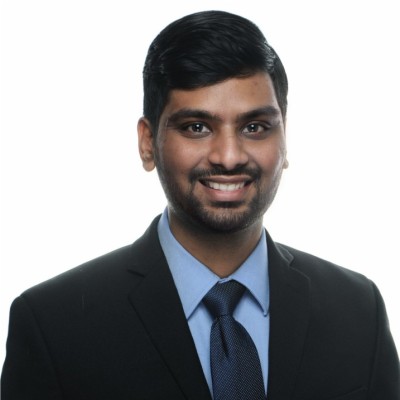}}]{Radhakrishna Tumbalam-Gooty}
received a B. Tech (Honors) and M. Tech degree in Chemical Engineering from Indian Institute of Technology Madras in Chennai, India. He received a Ph.D. degree in Chemical Engineering from Purdue University in West Lafayette, IN. He is now a Support Contractor for the National Energy Technology Laboratory, and a Senior Engineer at KeyLogic Systems in Pittsburgh, PA.
\end{IEEEbiography}

\begin{IEEEbiography}[{\includegraphics[width=1in,height=1.25in,clip,keepaspectratio]{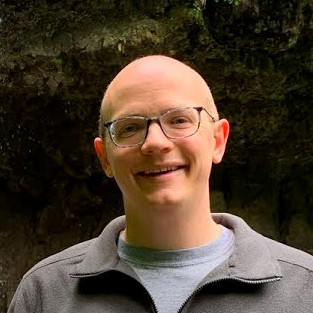}}]{Bernard Knueven}
received a B.S degree from Northern Kentucky University in Highland Heights, KY, and an M.S. degree from Miami University in Oxford, OH, both in Mathematics. He then received M.S. and Ph.D. degrees in Industrial Engineering from the University of Tennessee at Knoxville, TN. He is now a member of the Complex Systems Simulation \& Optimization group at the Computational Science Center at the National Renewable Energy Laboratory in Golden, CO. 
\end{IEEEbiography}

\begin{IEEEbiography}[{\includegraphics[width=1in,height=1.25in,clip,keepaspectratio]{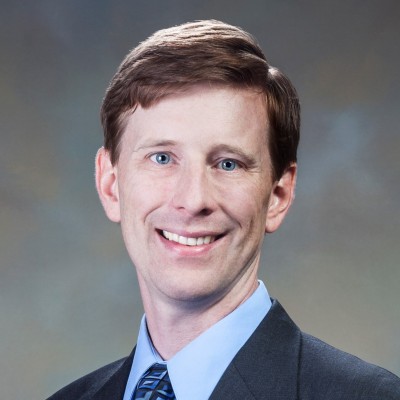}}]{John D. Siirola}
received a B.S. degree from Purdue University in West Lafayette, IN, and a Ph.D. degree from Carnegie Mellon University in Pittsburgh, PA, both in Chemical Engineering. He is now a member of the Technical Staff in the Discrete Math \& Optimization Department at Sandia National Laboratories in Albuquerque, NM. 
\end{IEEEbiography}

\begin{IEEEbiography}[{\includegraphics[width=1in,height=1.25in,clip,keepaspectratio]{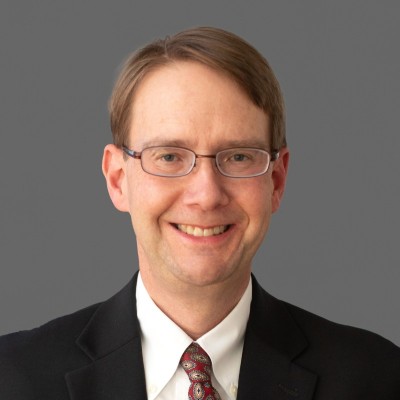}}]{David C. Miller}
received his B.S. degree from Rose-Hulman Institute of Technology in Terre Haute, IN, USA,  M.S. degree from the University of Illinois in Urbana-Champaign, IL, USA, and Ph.D. degree from The Ohio State University, Columbus, OH, USA. He spent 14 ½ years at the U.S. Department of Energy’s National Energy Technology Laboratory, serving in a variety of roles, including founding Technical Director of the Institute for the Design of Advanced Energy Systems (IDAES), Senior Fellow, and Chief Research Officer. He is currently the Chief Research Officer at OLI Systems, Inc. headquartered in NJ, USA.
\end{IEEEbiography}

\begin{IEEEbiography}[{\includegraphics[width=1in,height=1.25in,clip,keepaspectratio]{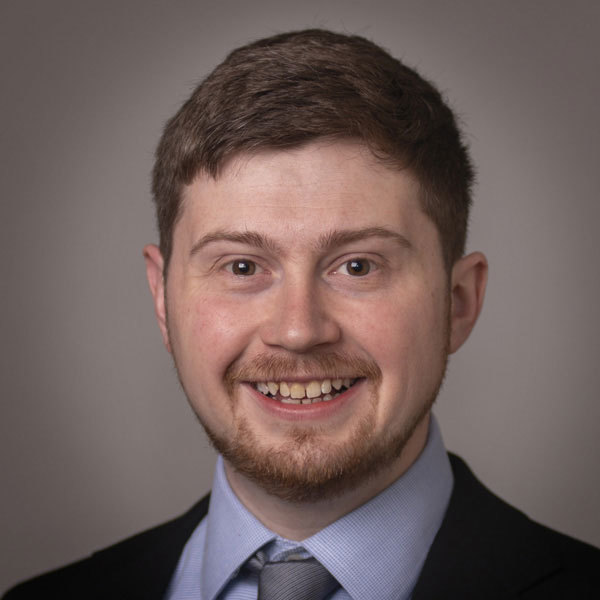}}]{Alexander W. Dowling}
received a B.S.E. degree from the University of Michigan, Ann Arbor, MI, and a Ph.D. degree from Carnegie Mellon University in Pittsburgh, PA, both in Chemical Engineering. He is now an associate Professor of Chemical and Biomolecular Engineering at the University of Notre Dame, Notre Dame, IN.
\end{IEEEbiography}


\vfill

\end{document}